# EDITOR'S COLUMN

### Richard S. Palais

The previous "Editors' Column", appearing in the January 1992 issue of the *Bulletin*, introduced a spirited debate on the merits of Information-based Complexity Theory. In that column we announced that the Research-Expository Papers might from time to time publish further debates on subjects of long-term interest to the general mathematical community. The rationale was that there should be at least one archival, scholarly journal of wide circulation willing to present discussions of controversial mathematical issues. To quote from what we said at that time:

> And not all such controversies are ancient history! Fifteen years ago there was a sharp controversy over purported excesses in the applications of Catastrophe Theory, and currently there is a similar controversy concerning what some see as an overselling and overpopularization of "fractals" and "chaos". Another simmering debate has grown out of the current renewal of the on-again, off-again love affair between Mathematics and Theoretical Physics. We have learned to accept that different standards of mathematical rigor may be appropriate when mathematics is being used as a tool to gain new insights about the physical world. But what standards should we apply to judge a paper that uses nonrigorous or semirigorous methods from physics to suggest important new insights into our own mathematical world, particularly if those insights seem beyond the reach of current rigorous mathematics?
>
> Especially because such questions cannot always be answered by logical principles alone, we believe that it is important for mathematicians to confront them. Even when rational discussion and debate does not completely resolve differences, at least it may clarify the issues.

In fact, that quotation prefigured the next such debate to appear in these pages. The July 1993 issue of the *Bulletin* contained an article by A. Jaffe and F. Quinn, with the title "Theoretical Mathematics: Towards a cultural synthesis of mathematics and theoretical physics", and responses to the Jaffe-Quinn article appear following these remarks. It is important for the reader to realize that these responses should not be read as some sort of referendum on the issues raised by Jaffe and Quinn. In the first place, those agreeing fully with an article write warm letters of congratulations to the authors, while those finding themselves in strong disagreement are more apt to write the editor asking for "equal time". In the present case, the editor's selection process favored dissent even more strongly than usual.

After the Jaffe-Quinn article was accepted, but before it actually appeared in print, it occurred to me that if I waited patiently for the article to be published and for people to respond to it on their own, then the time interval between the publication of the original article and the publication of the responses would be unacceptably long. Therefore, I wrote to various people, sending them an advance copy of the Jaffe-Quinn article and soliciting their comments. As for choosing whom to ask, since many people were mentioned by name in the Jaffe-Quinn article (some in ways that might be construed as critical), I felt it would be fair to ask all such if they would like to write a response. In addition, I wrote to others who, either by personal knowledge or hearsay, I realized had definite opinions on the matters addressed in the Jaffe-Quinn article. Finally, some of those whose opinions I solicited first wrote back suggesting the names of others I might wish to contact. I was personally highly gratified by the many thoughtful responses that came back, and I hope that the readers of the *Bulletin* will share in the pleasure and interest I took in reading them.

Starting with the January 1995 issue, the *Notices of the American Mathematical Society* will become much more archival in nature, and its Forum section will be an ideal place to carry on this kind of high-level debate on mathematical topics. Since this will remove the original rationale for publishing controversial articles in the *Bulletin*, I am (with mixed feelings) announcing the end of that program.